\documentclass[12pt]{article}
\usepackage{amsmath,amsthm,amsfonts,amssymb}

\newtheorem{theorem}{Theorem}
\newtheorem{lemma}{Lemma}

\newtheorem{prop}{Proposition}
\newtheorem{remark}{Remark}

\begin{document}

\title{A connection between the stochastic heat equation and fractional Brownian motion, and a simple proof of a result of Talagrand}

\author{Carl Mueller$^1$ and Zhixin Wu}
\date{}
\maketitle

\begin{abstract}
We give a new representation of fractional Brownian motion with Hurst 
parameter $H\leq\frac{1}{2}$ using stochastic 
partial differential equations.  This representation allows us to use the 
Markov property and time reversal, tools which are not usually available for 
fractional Brownian motion.  We then give simple proofs that fractional 
Brownian motion does not hit points in the critical dimension, and that it 
does not have double points in the critical dimension.  These facts were 
already known, but our proofs are quite simple and use some ideas of L\'evy.
\end{abstract}

\section{Introduction}
\label{intro}
\setcounter{equation}{0}

Our main result is a new representation for fractional Brownian motion using 
stochastic partial differential equations, described in this section.  As an 
application, in Section 2, Theorems \ref{th:hit} and \ref{th:double}, we state 
some known results about when fractional Brownian motion hits points and has 
double points.  Our representation allows us to give simple proofs of these 
results.  

In recent years there has been an upsurge of interest in fractional Brownian 
motion, see Nualart, Chapter 5 \cite{nua06}.  The most common model for noise 
in physical systems is white noise $\dot{B}_t$, the derivative of Brownian 
motion.  The central limit theorem gives some justification for using a 
Gaussian process such as $B_t$.  Furthermore, $\dot{B}_s,\dot{B}_t$ are 
independent if $s\ne t$.  In many situations, however, there are correlations 
between noise at different times.  A natural correlated Gaussian model to 
consider is fractional Brownian motion $X_t=X^H_t:t\ge0$ taking values in 
$\mathbf{R}^n$, with Hurst parameter $H\in(0,1]$.  The process $X_t$ is 
uniquely specified by the following axioms.
\begin{enumerate}
\item $X_0=0$ with probability 1.
\item $X_t:t\ge0$ is a Gaussian process with stationary increments.  That
is, for $t,h>0$ the probability distribution of the increment $X_{t+h}-X_t$
is independent of $t$.
\item For $c>0$ we have $X_{c t}\stackrel{\mathcal{D}}{=}c^HX_t$, where 
$\stackrel{\mathcal{D}}{=}$ denotes equality in distribution.  
\item $X_1$ has the standard normal distribution in $\mathbf{R}^n$.
\end{enumerate}
Note that Brownian motion is a fractional Brownian motion with Hurst parameter
$H=1/2$.

Next, we describe a seemingly unrelated process, the solution of the heat
equation with additive Gaussian noise.  Then we show that fractional Brownian
motion can be recovered from this solution.  There are several representations
of fractional Brownian motion, see Nualart \cite{nua06}, Chapter 5.  One
advantage of our representation is that we can use the Markov property and
time reversal, tools which fail for the fractional Brownian motion alone.
Using these extra tools, we give a simple proof of some hitting properties of 
fractional Brownian motion, and a result of Talagrand about double points.  
Throughout the paper we will write SPDE for ``stochastic partial differential 
equation''.

Let $N\geq1$.  Informally, we consider solutions 
$u(t,x):\mathbf{R}\times\mathbf{R}^N\to\mathbf{R}^n$ to the following 
equation.
\begin{eqnarray*}
\partial_tu &=& \Delta u + \dot{F}(t,x)  \\
u(-\infty,x) &=& 0
\end{eqnarray*}
where $\dot{F}(t,x)=(\dot{F}_1(t,x),\ldots,\dot{F}_n(t,x))$ is a generalized 
Gaussian field with the following covariance:
\[
E\left[\dot{F}_i(t,x)\cdot\dot{F}_j(s,y)\right] = \delta_{ij}\delta(t-s)h(x-y)
\]
where
\[
h(x) =
\left\{\begin{array}{ll}
\delta(x) & \mbox{if $H=\frac{1}{4}$} \\
|x|^{-\alpha} & \mbox{otherwise}
\end{array}
\right.
\]
and
\[
H = \frac{2-\alpha}{4}
\]
so that
\[
\alpha = 2-4H.  
\]
Furthermore, we have the following restrictions on $\alpha,H,N$.  
\begin{enumerate}
\item If $N=1$ then $0<\alpha\leq1$, so that $\frac{1}{4}\leq H<\frac{1}{2}$.  
\item If $N=2$ then $1<\alpha<2$, so that $0<H<\frac{1}{4}$.  	
\end{enumerate}
Note that if $H>\frac{1}{2}$ then $\alpha<0$ and $h(0)=0$, and then $h$ is
not a proper covariance.  
Our goal is to show that $X_t\stackrel{\mathcal{D}}{=}u(t,0)$, but this will 
not be literally true.  

The above description is not rigorous.  To be precise, $\dot{F}$ is a
centered Gaussian random linear functional on
$\mathbf{C}^\infty_c(\mathbf{R}^{N+1})$, the set of infinitely differentiable
functions with compact support on $(t,x)\in\mathbf{R}\times\mathbf{R}^{N}$,
taking values in $\mathbf{R}^n$, with covariance
\begin{equation}
\label{covar-fun-delta}
Q(f,g):= E\left[F(f)F(g)\right]
= \int_{\mathbf{R}}\int_{\mathbf{R}^N}f(t,x)\cdot g(t,x)dxdt
\end{equation}
if $H=\frac{1}{4}$, and
\begin{equation}
\label{covar-fun}
Q(f,g):=E\left[F(f)F(g)\right]
= \int_{\mathbf{R}}\int_{\mathbf{R}^N\times\mathbf{R}^N}
    f(t,x)\cdot g(t,y)h(x-y)dydxdt
\end{equation}
if $H\ne \frac{1}{4}$.  Note that in either case, 
the integral in (\ref{covar-fun-delta}) or (\ref{covar-fun}) is nonnegative 
definite.  Thus, we can extend $F(f)$ to all functions $f$ satisfying
\[
Q(f,f) < \infty.
\]
We call this class of functions $\mathbf{X}$.  Note that $\mathbf{X}$ 
implicitly depends on $\alpha,n,N$.  Furthermore, for $f$ taking values in 
$\mathbf{R}$, we say $f\in\mathbf{X}$ provided the $n$-dimensional vector 
$(f,\ldots,f)\in\mathbf{X}$.  

Next, for $t>0$ and $x\in\mathbf{R}^N$ let
\[
G(t,x) := 
\left\{\begin{array}{ll}
(4\pi t)^{-N/2}\exp\left(-\frac{|x|^2}{4t}\right)
        &\mbox{if $t>0$} \\
0 &\mbox{if $t\leq0$}
\end{array}
\right.
\]
be the heat kernel on $\mathbf{R}^N$.  

We would be tempted to define $u(t,x)$ by
\[
u(t,x) = \int_{-\infty}^{t}\int_{\mathbf{R}^N}G(t-s,x-y)F(dyds)
\]
but the integral will not converge.  However, $u(t,x)-u(0,0)$ looks more 
promising.  For $H=1/4$ and $N=1$, the stationary pinned string
was defined in \cite{mt01} as
\begin{equation}
\label{def-pinned}
U(t,x) :=
\int_{-\infty}^{t}\int_{\mathbf{R}}\left[G(t-s,x-y)-G(-s,-y)\right]F(dyds)
\end{equation}
when $t\geq 0$.  
This definition also works for other values of $H$, and $N\geq1$, provided
$g\in\mathbf{X}$, where $g(s,y)=g_{t,x}(s,y):=G(t-s,x-y)-G(-s,-y)$.

\begin{lemma}
\label{g-in-H}
Let $g$ be as in the previous paragraph.  For all $t\geq0$ and $x\in\mathbf{R}$ 
we have
\[
g(s,y)\mathbf{1}_{(s\leq t)} \in \mathbf{X}.
\]
\end{lemma}
We will prove Lemma \ref{g-in-H} in the Appendix.  

From the covariance of $\dot{F}$ one can easily deduce the following scaling 
property.  We leave the proof to the reader.  
\begin{lemma}
\label{scaling}
The noise $\dot{F}$ obeys the following scaling relation,
\[
\dot{F}(ct,c^{1/2}x)\stackrel{\mathcal{D}}{=} c^{-(2+\alpha)/4}\dot{F}(t,x).
\]
\end{lemma}

Turning to the SPDE, define $v(t,x)$ by $av(t,x)=U(ct,c^{1/2}x)$.  The reader 
can verify the following calculation using (\ref{def-pinned}).
\begin{eqnarray*}
a\partial_t v &=& c \partial_t U  \\
&=& c\left( \Delta U + \dot{F}(ct,c^{1/2}x) \right)   \\
&\stackrel{\mathcal{D}}{=}& a\Delta v + c^{1-(2+\alpha)/4}\dot{F}(t,x) \\
&\stackrel{\mathcal{D}}{=}& a\Delta v + c^{(2-\alpha)/4}\dot{F}(t,x)
\end{eqnarray*}
where the equality in distribution holds for the entire random field indexed 
by $t,x$.  Thus, we can cancel out the constants $c,a$ provided
\[
a=c^{\frac{2-\alpha}{4}}
\]
and then $v$ satisfies the same equation as $u$.  Thus,
\[
aU(t,x) \stackrel{\mathcal{D}}{=} U(ct,c^{1/2}x).
\]
Setting $x=0$ gives us the scaling relation for $U(t,0)$.  Thus we find
\begin{lemma}
\label{scaling-x}
$U(t,0)$ obeys the following scaling relation.  For $c>0$ we have
\[
U(ct,0) \stackrel{\mathcal{D}}{=} c^{\frac{2-\alpha}{4}}U(t,0)
\]
where the equality in distribution holds for the entire process indexed by $t$.
\end{lemma}

\begin{remark}
\label{rem:1}
Let $V_{t,x}(s,y)=U(t+s,x+y)-U(t,x)$.  It follows immediately from
(\ref{def-pinned}) that the random fields $V_{t,x}(s,y)$ and $U(s,y)$ are
equal in distribution.
\end{remark}

Let
\[
X_t = K_\alpha U(t,0)
\]
where 
\[
K_\alpha = \left[
\frac{(2-\alpha)\Gamma(\frac{n}{2})}
{2^{-\frac{3\alpha}{2}+1}\Gamma(\frac{n-\alpha}{2})}
\right]^{1/2}
\mbox{ if $\alpha\ne1$}
\]
and 
\[
K_\alpha = 2^{-1/2}(4\pi)^{d/4}
\mbox{ if $\alpha=1$.}
\]
We claim that
\begin{prop}
\label{prop:1}
Assume that $\alpha,N$ satisfy the conditions above, and let 
$X_t=K_\alpha U(t,0)$.  Then $X_t$, as defined above, is a fractional Brownian 
motion with Hurst parameter
\[
H = \frac{2-\alpha}{4}.
\]
\end{prop}

\noindent
\begin{proof}
We only need to verify the four axioms for fractional Brownian motion.
It follows from (\ref{def-pinned}) that $X_0=0$, so axiom 1 is 
satisfied.  Axiom 2 follows from Remark \ref{rem:1}.  Axiom 3 follows from the 
scaling properties of fractional Brownian motion and Lemma \ref{scaling-x}.  
Finally, Axiom 4 follows from (\ref{def-pinned}) and the integral of the 
covariance $h$, which we verify in the Appendix.  
\end{proof}

\textbf{Remark:}  Proposition \ref{prop:1} is related to a recent preprint
of Lei and Nualart \cite{ln08}.

\section{Critical dimension for hitting points, and for double points}
\setcounter{equation}{0}

The rest of the paper is devoted to the following questions.
\begin{enumerate}
\item For which values of $d,H$ does $X_t$ hit points?
\item For which values of $d,H$ does $X_t$ have double points?
\end{enumerate}
Recall that we say $X_t$ hits points if for each $z\in\mathbf{R}^n$, there
is a positive probability that $X_t=z$ for some $t>0$.  We say that $X_t$
has double points if there is a positive probability that $X_s=X_t$ for
some positive times $t\ne s$.  Here are our main results.  

\begin{theorem}
\label{th:hit}
Assume $0<H<\frac{1}{2}$, and that $\frac{1}{H}$ is an integer.  
For the critical dimension $n=\frac{1}{H}$, fractional Brownian motion does 
not hit points.  
\end{theorem}

\begin{theorem}
\label{th:double}
Assume $0<H<\frac{1}{2}$, and that $\frac{2}{H}$ is an integer.  
For the critical dimension $n=\frac{2}{H}$, fractional Brownian motion does 
not have double points.  
\end{theorem}

In fact, Talagrand answered the question of double points in \cite{tal98}, 
Theorem 1.1, and the assertion about hitting points was already known.  
Techniques from Gaussian processes, such as Theorem 22.1 of \cite{gh80}, can 
usually answer such questions except in the critical case, which is much more 
delicate.  The critical case is the set of parameters $n,H$ which lie on the 
boundary of the parameter set where the property occurs.  For example, 
$n=2,H=\frac{1}{2}$ falls in the the critical case for fractional Brownian 
motion to hit points.  But for $H=\frac{1}{2}$ we just have standard Brownian 
motion, which does not hit points in $\mathbf{R}^2$.  This illustrates the 
usual situation, that hitting does not occur, or double points do not occur, 
in the critical case.

It is not hard to guess the critical parameter set for fractional Brownian 
motion hitting points or having double points.  Heuristically, the range of a 
process with scaling $X_{\alpha t}\stackrel{\mathcal{D}}{=}\alpha^HX_t$ should 
have Hausdorff dimension $\frac{1}{H}$, if $X_t$ takes values in a space of 
dimension at least $\frac{1}{H}$.  For example, Brownian motion satisfies 
$B_{\alpha t}\stackrel{\mathcal{D}}{=}\alpha^{1/2}B_t$, and Brownian motion 
has range of Hausdorff dimension 2, at least if the Brownian motion takes 
values in $\mathbf{R}^n$ with $n\geq2$.  The critical parameter of $H$ for a 
process to hit points should be when the dimension of the range equals the 
dimension of the space.  Thus, the critical case for fractional Brownian 
motion taking values in $\mathbf{R}^n$ should be when the Hurst parameter is 
$H=1/n$.  For double points, we consider the 2-parameter process 
$V(s,t)=X_t-X_s$.  This process hits zero at double points of $X_t$, except 
when $t=s$.  The Hausdorff dimension of the range of $V$ should be 
$n=\frac{2}{H}$, and so the critical Hurst parameter for double points of $X_t$ 
should be $H=\frac{1}{2n}$.  

First note that the supercritical case can be reduced to the critical case.
That is, if $X_t=(X_t^{(1)}),\ldots,X_t^{(n+m)})=0$, then it is also true that
the projection $(X_t^{(1)}),\ldots,X_t^{(n)})=0$.  Furthermore, the
subcritical case is easier than the critical case, and it can be analyzed
using Theorem 22.1 of Geman and Horowitz \cite{gh80}.  Therefore, we
concentrate on the critical case $H=1/n$.

Below we give a simple argument inspired by \cite{mt01} which settles the 
critical case.  The argument goes back to L\'evy, and an excellent exposition is
given in Khoshnevisan \cite{kho03}.  It is based on scaling properties of the 
process, the Markov property, and time reversal.  Although fractional Brownian 
motion is not a Markov process, $U(t,x)$ does have the Markov property 
with respect to time.  Furthermore, it is time-reversible.

\section{Summary of L\'evy's argument}
\label{levy}
\setcounter{equation}{0}

Here is a brief summary of L\'evy's argument that 2-dimensional Brownian 
motion does not hit points.  Let $m(dx)$ denote Lebesgue measure on 
$\mathbf{R}^n$ and let $B_t$ denote Brownian motion on $\mathbf{R}^n$.  For 
this section, let $n=2$.  Furthermore, let $B[a,b]:=\{B_t: a\leq t\leq b\}$.  
It suffices to show that
\begin{equation}
\label{to-show}
E\Big[m\left(B[0,2]\right)\Big] = 0
\end{equation}
since then we would have
\begin{eqnarray*}
0 &=& E\left[\int_{\mathbf{R}^2}\mathbf{1}(z\in B[0,2])dz\right]   \\
&=& \int_{\mathbf{R}^2}P(z\in B[0,2])dz
\end{eqnarray*}
and so $P(z\in B[0,2])=0$ for almost every $z$.  

Next, for $0\leq t\leq 1$, let
\begin{eqnarray*}
Y_t &=& B_{1+t} - B_1  \\
Z_t &=& B_{1-t} - B_{1}.
\end{eqnarray*}
Recall that $Y_t,Z_t: 0\leq t\leq 1$ are independent standard 2-dimensional
Brownian motions.  This is a standard property of Brownian motion, which can 
be verfied by examining the covariances of $Y_t,Z_t: 0\leq t\leq 1$.  Then
$Y[0,1],Z[0,1]$ are independent random sets.  Furthermore, by Brownian scaling
and translation
\begin{eqnarray*}
E\Big[m\left(B[0,2]\right)\Big] 
&=& 2E\Big[m\left(B[0,1]\right)\Big] = 2E\Big[m\left(B[1,2]\right)\Big]  \\
&=& E\Big[m\left(B[0,1]\right)\Big] + E\Big[m\left(B[1,2]\right)\Big]  \\
&=&  E\Big[m\left(Y[0,1]\right)\Big] + E\Big[m\left(Z[0,1]\right)\Big].
\end{eqnarray*}
On the other hand, set theory gives us
\begin{eqnarray*}
\lefteqn{
E\Big[m\left(B[0,2]\right)\Big] = E\left[m(Y[0,1]\cup Z[0,1])\right]   }\\
&=& E\Big[m\left(Y[0,1]\right)\Big]
   + E\Big[m\left(Z[0,1]\right)\Big]
   - E\Big[m\left(Y[0,1]\cap Z[0,1]\right)\Big]
\end{eqnarray*}
and therefore
\[
E\Big[m\left(Y[0,1]\cap Z[0,1]\right)\Big] = 0.
\]
By Fubini's theorem,
\begin{eqnarray}
\label{to-show-2}
0 &=& E\Big[m\left(Y[0,1]\cap Z[0,1]\right)\Big]   \\
&=& E\left[\int_{\mathbf{R}^2} \mathbf{1}_{Y[0,1]}(z)
       \mathbf{1}_{Z[0,1]}(z)dz\right]  \nonumber\\
&=& \int_{\mathbf{R}^2} E\left[\mathbf{1}_{Y[0,1]}(z)
       \mathbf{1}_{Z[0,1]}(z)\right]dz.    \nonumber
\end{eqnarray}
By the independence of $Y[0,1]$ and $Z[0,1]$ and the Cauchy-Schwarz 
inequality, we have
\begin{eqnarray}
\label{to-show-3}
0 &=& \int_{\mathbf{R}^2} E\left[\mathbf{1}_{Y[0,1]}(z)
       \mathbf{1}_{Z[0,1]}(z)\right]dz  \\
&=& \int_{\mathbf{R}^2} E\left[\mathbf{1}_{Y[0,1]}(z)\right]
       E\left[\mathbf{1}_{Z[0,1]}(z)\right]dz  \nonumber\\
&=& \int_{\mathbf{R}^2} \left(E\left[\mathbf{1}_{Y[0,1]}(z)\right]\right)^2
  dz  \nonumber\\
&\geq& \left(\int_{\mathbf{R}^2} E\left[\mathbf{1}_{Y[0,1]}(z)\right]
  dz\right)^2  \nonumber\\
&=& \Big(Em(Y[0,1])\Big)^2.  \nonumber
\end{eqnarray}
Therefore $E[m(Y[0,1])]=0$ and (\ref{to-show}) follows from the definition of 
$Y$.

\section{Proof of Theorems \ref{th:hit} and \ref{th:double}}
\setcounter{equation}{0}

Now we use L\'evy's argument to prove our main theorems.  Throughout, we 
assume that $H,\alpha,N$ satisfy the restrictions given in the introduction.  

\subsection{Hitting points, Theorem \ref{th:hit}}

The argument exactly follows that in Section \ref{levy}, except that
$\mathbf{R}^2$ is replaced by $\mathbf{R}^n$.  Also, by axiom (3),
\[
X_{ct} \stackrel{\mathcal{D}}{=} c^{H} X_t
\]
for $c>0$.  However, since $X_t$ takes values in $\mathbf{R}^n$ and $H=1/n$, 
we still have
\[
E\Big[m\left(X[0,2]\right)\Big] = 2 E\Big[m\left(X[0,1]\right)\Big]
= 2 E\Big[m\left(X[1,2]\right)\Big].
\]
Recall that $m(\cdot)$ denotes Lebesgue measure in $\mathbf{R}^n$.  
As before, let
\begin{eqnarray*}
Y_t &=& X_{1+t} - X_1  \\
Z_t &=& X_{1-t} - X_{1}.
\end{eqnarray*}
It is no longer true that $Y[0,1],Z[0,1]$ are independent.  Now we
use the fact that $X_t$ is equal in distribution to $u(t,0)$,
where $u(t,x)$ is the stationary pinned string.  Changing the
probability space if necessary, let us write $X_t=K_\alpha U(t,0)$, and let
$\mathcal{H}_t$ denote the $\sigma$-field generated by $U(t,x):
x\in\mathbf{R}^n$.  Then we have the following lemma.  

\begin{lemma}
\label{thesis-1}
Let us use the above notation.  Then $Y[0,1],Z[0,1]$ are conditionally 
independent and identically distributed given $\mathcal{H}_1$.  
\end{lemma}

\begin{proof}[Proof of Lemma \ref{thesis-1}.]
The lemma is proved in \cite{mt01}, Corollary 1, for $H=1/4$, and the proof 
for other values of $H$ uses similar ideas.  

To show that $Y[0,1],Z[0,1]$ are identically distributed, we merely use the 
definition of $Y,Z$ and make a change of variable.  Below, equality in 
distribution means that the processes indexed by $t$ are equal in 
distribution. 
\begin{eqnarray*}
Y_t &=& X_{1+t} - X_1  \\
&=& K_\alpha U(1+t,0) - K_\alpha U(1,0)  \\
&=& 
K_\alpha\int_{-\infty}^{1+t}\int_{\mathbf{R}^N}
    \left[G(1+t-s,-y)-G(1-s,-y)\right]F(dyds)   \\
&\stackrel{\mathcal{D}}{=}&
K_\alpha\int_{-\infty}^{t}\int_{\mathbf{R}^N}
    \left[G(t-s,-y)-G(-s,-y)\right]F(dyds).
\end{eqnarray*}
But by the definition of $Z_t$,
\begin{eqnarray*}
Z_t &=& X_{1-t} - X_{1}  \\
&=& K_\alpha \Big(U(1-t,0) - U(1,0)\Big)    \\
&=& 
-K_\alpha\int_{-\infty}^{1}\int_{\mathbf{R}^N}
    \left[G(1-s,-y)-G(1-t-s,-y)\right]F(dyds)   \\
&\stackrel{\mathcal{D}}{=}&
K_\alpha\int_{-\infty}^{t}\int_{\mathbf{R}^N}
    \left[G(t-s,-y)-G(-s,-y)\right]F(dyds)
\end{eqnarray*}
and therefore $Y_t,Z_t$ are identically distributed processes.  

Next we discuss the conditional independence of $Y[0,1],Z[0,1]$ given 
$\mathcal{H}_1$.  First we claim that the stationary pinned string $U(t,x)$ 
enjoys the Markov property with respect to $t$.  This is a general fact about 
stochastic evolution equations, and we refer the reader to \cite{wal86}, 
Chapter 3.  It follows that $Z[0,1]$ is conditionally independent of 
$Y[0,1]$ given $\mathcal{H}_1$.  

This proves Lemma \ref{thesis-1}.
\end{proof}

From here we duplicate the argument in Section \ref{levy}, replacing
expectation by conditional expectation given $\mathcal{H}_1$.  Briefly, it
suffices to show that
\begin{equation}
\label{to-show-1}
E\Big[m\left(Y[0,1]\right)\Big|\mathcal{H}_1\Big] = 0.
\end{equation}
But, following the same argument as before, we conclude that with probability
one,
\[
E\Big[m\left(Y[0,1]\cap Z[0,1]\right)\Big|\mathcal{H}_1\Big] = 0.
\]
We leave it to the reader to verify that (\ref{to-show-2}) and
(\ref{to-show-3}) still hold, provided expectation is replaced by conditional
expectation given $\mathcal{H}_1$.  This verifies (\ref{to-show-1}), and 
finishes the proof of Theorem \ref{th:hit}.

\subsection{Double points, Theorem \ref{th:double}}

To show that $X_t$ does not have double points in the critical case, we use
the same argument, but applied to the two-parameter process
\[
V(s,t) := X(t) - X(s).
\]
We need to show that $V(s,t)$ has no zeros except if $s=t$.  To simplify the
argument, we will show that $V(s,t)$ has no zeros for $(s,t)\in\mathcal{R}$,
where
\[
\mathcal{R} := [0,2]\times[4,6].
\]
The same argument would apply to any other rectangle whose intersection with
the diagonal has measure 0.  Let us subdivide $\mathcal{R}$ into 4
subrectangles $\mathcal{R}_i: i=1,\ldots,4$ each of which is a translation of
$[0,1]^2$.  Again we argue as in Section \ref{levy}.  By scaling, we see that
for each $i=1,\ldots,4$
\[
E\left[m(V(\mathcal{R}))\right] = 4E\left[m(V(\mathcal{R}_i))\right].
\]
Next, let $\mathcal{H}_1$ be the $\sigma$-field generated by
$\{u(1,x): x\in\mathbf{R}^n\}$, and suppose we have labeled the $\mathcal{R}_i$
such that $\mathcal{R}_1=[0,1]\times[4,5]$ and $\mathcal{R}_2=[1,2]\times[4,5]$.
Thus, as before, for each pair $i\ne j\in\{1,\ldots,4\}$ we have
\[
E\Big[m(V(\mathcal{R}_i)\cap V(\mathcal{R}_j))\Big|\mathcal{H}_1\Big] = 0.
\]
Now in \cite{mt01}, Corollary 1, it was shown that for $H=1/4$,
$V(\mathcal{R}_1),V(\mathcal{R}_2)$ are  conditionally i.i.d. given
$\mathcal{H}_1$.  For other values of $H$, the argument is very similar to the 
proof of Lemma \ref{thesis-1}, and we leave the details to the reader.  

Therefore, as in Section \ref{levy}, we conclude that with
probability one,
\[
E\Big[m(V(\mathcal{R}_i))\Big|\mathcal{H}_1\Big] = 0.
\]
Also as in Section \ref{levy}, this finishes the proof of Theorem 
\ref{th:double}.

\appendix

\section{Appendix}
\setcounter{equation}{0}

We first recall a standard fact about the Fourier transform in $\mathbf{R}^N$, 
which we take from Lemma 4.1 of Wolff \cite{wol03}.
\begin{lemma}
\label{l:fourier}
Let $h_a(x)=\frac{\gamma(a/2)}{\pi^{a/2}}|x^{-a}|$.  Then $\hat{h_a}=h_{N-a}$ 
in the sense of $\mathbf{L}^1+\mathbf{L}^2$ Fourier transforms if 
$\frac{N}{2}<\mbox{Re}(a)<N$, and in the sense of distributional Fourier 
transforms if $0<\mbox{Re}(a)<N$.  Here $\gamma$ is the gamma function.  
\end{lemma}

Taking the Fourier transform in the distributional sense is enough, because we 
can use cutoffs and then take limits.  

Now we give the proof of Lemma \ref{g-in-H}.

\begin{proof}[Proof of Lemma \ref{g-in-H}.]
In the case of $H=\frac{1}{4}$, for which $h(x)=\delta(x)$,
Lemma \ref{g-in-H} follows from Proposition 1 of \cite{mt01}.

Next we move on to the case of $H\ne\frac{1}{4}$.  By using the triangle 
inequality and changing variables, and scaling, we see that it suffices to 
prove the following inequalities for all $t>0$ and $x\in\mathbf{R}$.  
\begin{align}
\label{eq:1}
\int_{0}^{\infty}\int_{\mathbf{R}^N\times\mathbf{R}^N}&
  \Big(G(s,z+1)-G(s,z)\Big)^2   \\
&  \times\Big(G(s,z'+1)-G(s,z')\Big)^2
  |z-z'|^{-\alpha}dzdz'ds < \infty  \nonumber\\
\label{eq:2}
\int_{0}^{\infty}\int_{\mathbf{R}^N\times\mathbf{R}^N}&
  \Big(G(s+1,z)-G(s,z)\Big)^2|z-z'|^{-\alpha}   \\
&  \times\Big(G(s+1,z')-G(s,z')\Big)^2
    |z-z'|^{-\alpha}dzdz'ds < \infty  \nonumber\\
\label{eq:3}
\int_{0}^{1}\int_{\mathbf{R}^N\times\mathbf{R}^N}&
  G(s,z)G(s,z')|z-z'|^{-\alpha}dzdz'ds < \infty 
\end{align}

First we deal with (\ref{eq:2}).  Taking the Fourier transform of $G(s,z)$ 
with respect to $z$, we recall that $\hat{G}(s,\xi)=\exp(-s|\xi|^2)$.  Also by 
Lemma \ref{l:fourier} and our restrictions on $\alpha,N$, the Fourier 
transform of $|x|^{-\alpha}$ is $c|\xi|^{\alpha-N}$ for some finite constant 
$c$.  Then Plancherel's theorem and Fubini's theorem show that (\ref{eq:2}) 
equals a constant times
\begin{align*}
\int_{0}^{\infty}\int_{\mathbf{R}^N} 
&\left(e^{-(s+1)|\xi|^2}-e^{-s|\xi|^2}\right)^2
  |\xi|^{\alpha-N}d\xi ds   \\
&= \int_{\mathbf{R}^N}\int_{0}^{\infty}
  e^{-2s|\xi|^2}\left(e^{-|\xi|^2}-1\right)^2
  |\xi|^{\alpha-N}dsd\xi    \\
&= \frac{1}{2}\int_{\mathbf{R}^N}
   \left(e^{-|\xi|^2}-1\right)^2
  |\xi|^{\alpha-N-2}d\xi   \\
&<\infty.  
\end{align*}
The final inequality can be verified by noting the restrictions on $\alpha,N$, 
splitting up the preceding integral into integrals over $|\xi|<1$ and 
$|\xi|\geq1$, using the bound $|e^{-\xi^2}-1|\leq \min(1,\xi^2)$, and 
switching to polar coordinates.

Secondly we treat with (\ref{eq:3}).  Using the Fourier transform as in the 
previous case, we find that (\ref{eq:3}) equals a constant times
\begin{equation*}
\int_{\mathbf{R}^N}e^{-2s|\xi|^2}|\xi|^{\alpha-N}d\xi
= Cs^{N-\alpha-1}
\end{equation*}
In view of our restrictions on $\alpha,N$, we see that
\begin{equation*}
\int_{0}^{1}s^{N-\alpha-1}ds < \infty
\end{equation*}
and so (\ref{eq:3}) is finite.  

We can set $z=ys^{1/2}$ with $s$ fixed to 
deduce
\begin{equation}
\label{detail:1}
\int_{\mathbf{R}^N\times\mathbf{R}^N}G(s,z)G(s,z')|z-z'|^{-\alpha}dzdz'
= Cs^{-\alpha/2}
\end{equation}
where the reader can check that $C<\infty$.  Because of our restrictions on 
$\alpha,N$, we see that the integral (\ref{detail:1}) over $s\in[0,1]$ is 
finite, verifying (\ref{eq:3}).  

Finally, we treat (\ref{eq:1}).  We use the preceding facts about the Fourier 
transform, and also the fact that the Fourier transform of $f(x+a)$ is 
$\hat{f}(\xi)e^{ia\xi}$.  By Plancherel's theorem, we find that 
that (\ref{eq:1}) equals a constant times
\begin{equation*}
\int_{0}^{\infty}\int_{\mathbf{R}^N} 
  e^{-2s|\xi|^2}\left(e^{-|\xi|^2}-1\right)^2
  |\xi|^{\alpha-1}d\xi ds
<\infty
\end{equation*}
by the same reasoning as for (\ref{eq:2}).

\end{proof}

\end{document}